\documentclass[12pt,a4paper]{article}
\usepackage{graphicx}
\usepackage{graphics}
\usepackage[bf]{caption2}
\usepackage{graphicx}
\usepackage{graphics}
\usepackage{amssymb}
\usepackage{amsmath}
\usepackage[cp1251]{inputenc}

\newcommand{\vp}{{\mathbf p}}

\newcommand{\vx}{{\mathbf x}}

\newcounter{theoremcounter}

\newcounter{remarkcounter}

\newcounter{corollarycounter}
\newtheorem{theorem}[theoremcounter]{Theorem}

\newtheorem{remark}[remarkcounter]{Remark}
\newtheorem{corollary}[corollarycounter]{Corollary}

\begin{document}


\title{On the limiting characteristics for an inhomogeneous $M_t|M_t|S$ queue with catastrophes} 

\date{}

\author{Alexander Zeifman    \and Anna Korotysheva \and  \  \and Victor Korolev }

\author{{ \small Alexander Zeifman,   Vologda State  University,}\\ {\small
Institute of Informatics Problems RAS, and ISEDT RAS}
\\ { \small Anna Korotysheva,     Vologda State University,}
\\ { \small Victor Korolev,  Moscow State University}\\ {\small and Institute of Informatics Problems RAS}
}
\maketitle 

{\bf Abstract.} We study weak ergodicity, bounds on the rate of
convergence,  and problems of computing of the limiting
characteristics for an inhomogeneous $M_t|M_t|S$ queueing model with
possible catastrophes.



\section{Introduction}

Qualitative and quantitative  properties of inhomogeneous
continuous-time Markov chains and the correspondent queueing models
have been investigated since 1980's, see for instance first results
in \cite{z85},  \cite{z94}, \cite{z95}. Queueing systems with
catastrophes (queues with disasters) in different situations were
studied by a number of authors, see, for instance,
\cite{Di08},\cite{du1,du2,du3}, \cite{z09a,z09b}.

Perturbation bounds for an inhomogeneous $M_t|M_t|N$ queue with
catastrophes were obtained on \cite{z12}. First investigations for
this model with catastrophes rates depending on the length of the
queue were studied in \cite{z09b}.

It is well known that explicit expressions for the probability
characteristics of stochastic models can be found only in a few
special cases. If we deal with inhomogeneous Markovian model, then
we must approximately calculate, in addition, the limiting
probability characteristics of the process. The problem of existence
and construction of limiting characteristics for time-inhomogeneous
birth and death processes is important for queueing applications,
see for instance,
\cite{d95,d03,d12,gz04,man95,mas94,mas13,tan13,z06}. A general
approach to the study of the rate of convergence for birth-death
models and related bounds were considered in \cite{z95}, and for
finite birth-death-catastrophe models they were considered in
\cite{z13a}. Calculation of the limiting characteristics for the
process via truncations was firstly mentioned in \cite{z88} and was
considered in details in \cite{z06}. The best results in this
direction for general inhomogeneous birth-death models were obtained
in our recent paper \cite{zsks14}.

Here we apply this general approach to an inhomogeneous $M_t|M_t|S$
queue with catastrophes in a general situation where catastrophe
rates depend on the length of the queue. Moreover, we will obtain
and discuss explicit bounds on the rate of convergence to the
limiting characteristics in weak ergodic situation as well as
approximation bounds of the limiting characteristics. Finally, we
discuss an example of this queueing model.

Let $X=X(t)$, $t\geq 0$, be an inhomogeneous, in general,
continuous-time Markov chain, which is the queue length process for
the corresponding queueing model.

Let $p_{ij}(s,t)=Pr\left\{ X(t)=j\left| X(s)=i\right. \right\}$,
$i,j \ge 0, \;0\leq s\leq t$ be the transition probabilities for
$X=X(t)$,  $p_i(t)=Pr\left\{ X(t) =i \right\}$  be its state
probabilities, and ${\bf p}(t) = \left(p_0(t), p_1(t),
\dots\right)^T$ be the corresponding probability distribution.

In the inhomogeneous case we assume that all intensity functions are
linear combinations of a finite number of  nonnegative functions
which are locally integrable on $[0,\infty)$. Then the corresponding
transposed intensity matrix is
$$
A(t)=\left(
\begin{array}{cccccccc}
a_{00}(t) & \mu_1(t)+\xi_1(t) & \xi_2(t)  & \xi_3(t)  & \xi_4(t) & \xi_5(t)  & \ldots\\[3pt]
\lambda_0(t)    & a_{11}(t)  & \mu_{2}(t)  & 0   & 0  & 0 & \ldots\\[3pt]
0 & \lambda_1(t) & a_{22}(t)&  \mu_{3}(t)  & 0    &  0 &   \ldots\\[3pt]
0 & 0 & \lambda_2(t) & a_{33}(t)&  \mu_{4}(t)  & 0    &   \ldots\\[3pt]
\ldots&\ldots&\ldots&\ldots&\ldots&\ldots&\ldots\\[3pt]
\end{array}
\right),
$$
\noindent  $a_{ii}(t) =-\sum_{j\neq i} a_{ji}(t)$. In addition,
applying our standard approach (see details in
\cite{gz04,z95,z08,z06}) we suppose that the intensity matrix is
essentially bounded, i. e.
\begin{equation}
|a_{ii}(t)| \le L < \infty,
\label{0102-1}
\end{equation}
\noindent for almost all $t \ge 0$. Then the probabilistic dynamics
of the process is represented by the forward Kolmogorov system
\begin{equation}
 \left\{
\begin{array}{cc}
\frac{dp_0}{dt} = -\lambda_0 (t) p_0 +\mu_1(t)p_1 + \sum_{k \ge 1}\xi_k(t)p_k ,  \\
\frac{dp_k}{dt} = \lambda_{k-1} (t)p_{k-1} -\left(\lambda_{k} (t) +
 \mu_{k}(t)+
\xi_k(t)\right)p_k +\mu_{k+1}(t) p_{k+1},  k \ge 1,
\end{array}
\right.  \label{1001}
\end{equation}
\noindent where $\lambda_{k}(t) = \lambda(t)$,  $\mu_{k}(t) =
\min(k,S)\mu(t)$, and $\xi_k(t)$ are the arrival, service and
catastrophe rates, respectively.

Throughout the paper by $\|\,\cdot\,\|$  we denote  the $l_1$-norm,
i. e.,  $\|{\vx}\|=\sum|x_i|$, and $\|B\| = \sup_j \sum_i |b_{ij}|$
for $B = (b_{ij})_{i,j=0}^{\infty}$. Let $\Omega$ be the set all
stochastic vectors, i. e. $l_1$-vectors with nonnegative coordinates
and unit norm.

Then we have $\|A(t)\| = 2\sup_{k}\left|a_{kk}(t)\right| \le 2 L $
for almost all $t \ge 0$. Hence, the operator function $A(t)$ from
$l_1$ into itself is bounded for almost all $t \ge 0$ and locally
integrable on $[0;\infty)$. Therefore, we can consider the forward
Kolmogorov system (\ref{1001})
\begin{equation} \label{ur01}
\frac{d\vp}{dt}=A(t)\vp(t),
\end{equation}
as a differential equation in the space $l_1$ with bounded operator.

It is well known (see \cite{DK}) that the Cauchy problem for
differential equation (\ref{ur01}) has a unique solution for an
arbitrary initial condition, and  $\vp(s) \in \Omega$ implies
$\vp(t) \in \Omega$ for $t \ge s \ge 0$.

By $E(t,k) = E\left\{X(t)\left|X(0)=k\right.\right\}$ denote the
mean (the mathematical expectation) of the queue length process
$X(t)$ at the moment $t$ under the initial condition $X(0)=k$.

Recall that a Markov chain $X(t)$ is called {\it weakly ergodic}, if
$\|{\bf p}^*(t)-{\bf p}^{**}(t)\| \to 0$ as $t \to \infty$ for any
initial conditions ${\bf p}^*(0), {\bf p}^{**}(0)$, where ${\bf
p}^*(t)$ and ${\bf p}^{**}(t)$ are the corresponding solutions of
(\ref{ur01}). A Markov chain $X(t)$ has  the {\it limiting mean}
$\varphi (t)$, if $ \lim_{t \to \infty }  \left(\varphi (t) -
E(t,k)\right) = 0$ for any $k$.

\section{Ergodicity bounds}

Consider an increasing sequence of positive numbers  $\{d_i\}$,
$i=0,1,2, \dots$, $d_0=1$, and the corresponding triangular matrix
$D$:
\begin{equation}
D=\left(
\begin{array}{ccccccc}
d_1   & d_1 & d_1 & \cdots  \\
0   & d_2  & d_2  &   \cdots  \\
0   & 0  & d_3  &   \cdots  \\
& \ddots & \ddots & \ddots \\
\end{array}
\right) \label{204}
\end{equation}
Let $l_{1D}$ be the space of sequences:
$$
l_{1D}=\left\{{\bf z} = (p_1,p_2,\cdots)^T :\, \|{\bf z}\|_{1D}
\equiv \|D {\bf z}\| <\infty \right\}. $$
We also introduce the
auxiliary space of sequences $l_{1E}$ as
$$l_{1E}=\left\{{\bf z} = (p_1,p_2,\cdots)^T :\, \|{\bf z}\|_{1E}
\equiv \sum k|p_k| <\infty \right\}.$$

Put $$W=\inf_{i \ge 1} \frac {d_i}{i}, \quad g_i=\sum_{n=1}^i d_n.$$
Consider the following expressions:
\begin{eqnarray}
\alpha_{k}\left( t\right) = \lambda _k\left( t\right) +\mu
_{k+1}\left( t\right) +\xi_{k+1}(t) - \frac{d_{k+1}}{d_k} \lambda
_{k+1}\left( t\right) - \nonumber \\ \frac{d_{k-1}}{d_k} \mu
_k\left( t\right), \quad k \ge 0, \label{211}
\end{eqnarray}
and
\begin{equation}
\alpha\left( t\right) = \inf_{k\geq 0} \alpha_{k}\left( t\right).
\label{212}
\end{equation}

Now recall the following general statement.

\begin{theorem}\label{t1}
Let $X(t)$ be a birth-death-catastrophe process (BDPC)  with rates
$\lambda_k(t)$, $\mu_k(t)$ and $\xi_{k}(t)$. Assume that there
exists a sequence $\{d_i\}$ such that
\begin{equation}
\int\limits_0^\infty \alpha(t)\, dt = + \infty . \label{213}
\end{equation}
Then $X(t)$ is weakly ergodic, and the following bounds hold:
\begin{equation}
\|{\bf p^*}(t) - {\bf p^{**}}(t)\|_{1D}  \le e^{-\int\limits_s^t
\alpha(\tau)\, d\tau}\|{\bf p^*}(s) - {\bf p^{**}}(s)\|_{1D},
\label{214}
\end{equation}
\begin{equation}
\|{\bf p^*}(t) - {\bf p^{**}}(t)\|  \le 4 e^{-\int\limits_s^t
\alpha(\tau)\, d\tau}\sum_{i \ge 1}g_i |p^*_i(s) - p^{**}_i(s)|,
\label{215}
\end{equation}
\noindent for any $ t \ge s \ge 0$ and any initial conditions ${\bf
p^*}(s), {\bf p^{**}}(s)$. \label{th01}
\end{theorem}

\textbf{Proof.} The proof follows the lines of the reasoning used to
prove Theorem 3 in \cite{zsks14}, hence we only outline this
argumentation here. Put $p_0(t) = 1 - \sum_{i \ge 1} p_i(t)$, then
from (\ref{ur01}) we have the following system:
\begin{equation}\label{216}
\frac{d{\bf z}(t)}{dt}=B(t){\bf z}(t)+{\bf f}(t),
\end{equation}
\noindent  where ${\bf z(t)} = \left(p_1(t),p_2(t),\dots\right)^T,$
${\bf f(t)} = \left(\lambda_0(t),0,0,\dots\right)^T$, $B(t) =
\left(b_{ij}(t)\right)_{i,j=1}^{\infty}$ and
\begin{equation}
b_{ij} = \left\{
\begin{array}{rlccccc}
 -(\lambda_0+\lambda_1+\mu_1+\xi_{1}),  & \mbox { if } \quad i=j=1, \\
\mu_2-\lambda_0,  & \mbox { if } \quad i=1,j=2, \\
-\lambda_0,  & \mbox { if } \quad i=1,j>2, \\
-(\lambda_j+\mu_j+\xi_{j}),  & \mbox { if } \quad i=j >1, \\
\mu_j,  & \mbox { if } \quad i=j-1 >1, \\
\lambda_j,  & \mbox { if } \quad i=j+1 >1, \\
0,  & \mbox { otherwise. }
\end{array}
\right.  \label{217}
\end{equation}
This is a linear non-homogeneous differential system, the solution
of which can be written as
\begin{eqnarray}
{\bf z}(t)&=& V(t,0){\bf z}(0)+\int_0^t{}V(t,\tau){\bf
f}(\tau)\,d\tau, \label{218}
\end{eqnarray}
where $V(t,z)$ is the Cauchy operator of (\ref{216}), see, for
instance,  \cite{zsks14}.

We can consider (\ref{216}) as a differential equation in the space
$l_{1D}$ with bounded and locally integrable on $[0,\infty)$
coefficients  ${\bf f}(t)$ and $B(t)$.

Applying the notion of the logarithmic norm and the related bounds
(see  \cite{dzp,gz04,z06,zsks14} for details), we obtain  the
following bound for  the logarithmic norm $\gamma \left(B(t)\right)$
in $l_{1D}$: {
\begin{eqnarray}
\gamma \left(B\right)_{1D} =  \gamma \left(DB(t)D^{-1}\right)_1 =  \nonumber \\
\sup\limits_{i \ge 0}\left(\frac{d_{i+1}}{d_i}
 \lambda_{i+1}(t) + \frac{d_{i-1}}{d_i} \mu_i(t) - \left(\lambda_i(t)
 +\mu_{i+1}(t)+\xi_{i+1}(t)
\right)   \right) \nonumber  =
\\  - \inf \limits_{k \ge 0} \left(\alpha_{k}\left( t\right)\right)  = -
\alpha(t),  \label{219}
\end{eqnarray}}
\noindent in accordance with (\ref{212}). Hence,
\begin{equation}
\|V(t,s)\|_{1D} \le e^{-\int\limits_s^t \alpha(\tau)\, d\tau}.
\label{220}
\end{equation}
Therefore, bound (\ref{214}) takes place.

On the other hand, inequalities $ \|{\bf z}\| \le 2\|{\bf
z}\|_{1D}$,  and $\|{\bf p^*} - {\bf p^{**}}\| \le 2\|{\bf z}\|$ for
any ${\bf p^*},  {\bf p^{**}}$ and corresponding ${\bf z}$ (see, for
instance \cite{zsks14}) imply bound (\ref{215}).

\smallskip

\begin{corollary}{}
Let, in addition, the numbers $d_i$ grow sufficiently fast so that
$W > 0$. Then $X(t)$ has the limiting mean, say $\phi(t)$, and the
following bound holds:
\begin{equation}
\left|\phi(t) - E_k(t)\right| \le \frac{4}{W} e^{-\int\limits_0^t
\alpha(\tau)\, d\tau}\|{\bf p}(0) - {\bf e_k}\|_{1D}. \label{222}
\end{equation}
\end{corollary}

\smallskip

Now we can obtain ergodicity bounds for the queue length process of
an $M_t|M_t|S$ queue with catastrophes.

\begin{theorem} \label{t2}
Let $\xi_k(t)=\zeta_k \xi(t)$,
\begin{equation}
\inf_n \zeta_n = \zeta > 0 \label{zkts02}
\end{equation}
\noindent and let there exist $\varepsilon > 0$ such that
\begin{equation}
\int_0^{\infty} \left(\zeta \xi(t)-\varepsilon\lambda(t)\right)\, dt
= +\infty \label{zkts02a}
\end{equation}
$($large catastrophe rates$)$. Then the queue-length process $X (t)
$ is weakly ergodic, has the limiting mean, and the following bounds
hold:
\begin{equation}
\|{\bf p}(t) - {\bf \pi}(t)\| \le 4\left(1+\varepsilon\right)^k \varepsilon^{-1} e^{-\int\limits_0^t
\left(\zeta\xi(\tau) - \varepsilon \lambda(\tau)\right)\, d\tau} \label{zkts0201}
\end{equation}
\noindent
\begin{equation}
\left|E_k(t) - E_0(t)\right| \le
\frac{4\left(1+\varepsilon\right)^k}{\varepsilon W}
e^{-\int\limits_0^t \left(\zeta\xi(\tau) - \varepsilon
\lambda(\tau)\right)\, d\tau}, \label{zkts0202}
\end{equation}
\noindent for any initial number of customers $ X (0) = k $, where $
{\bf \pi} (t) $ and $ \phi (t)= E_0(t) $ are the limiting regime and
the limiting mean correspondent to the empty initial length of the
queue.
\end{theorem}
\smallskip

{\bf Proof.} Put $d_0=1$, $d_{k+1}=(1+\varepsilon)d_k,\ k \ge 0$,
then instead of (\ref{219}) we have the  following bound for the
logarithmic norm:
\begin{eqnarray}
\gamma \left(B(t)\right)_{1D} \le  -\left(\zeta\xi(t) + \frac{\varepsilon}{1+\varepsilon}\left(S\mu(t)-(1+\varepsilon)\lambda(t)\right)\right) \le \\\nonumber
 -\left(\zeta\xi(t)  - \varepsilon\lambda(t)\right)  = -\alpha_*(t), \label{zkts040}
\end{eqnarray}
\noindent where  $\int_0^{\infty} \alpha_*(t) \, dt = + \infty$ in
accordance with (\ref{zkts02}).

Putting ${\bf p}^*(0) = {\bf \pi}(0)= {\bf e}_0$, ${\bf p}^{**}(0) =
{\bf p}(0)= {\bf e}_k$, from Theorem 1  we obtain bound
(\ref{zkts0201}). The second estimate follows from (\ref{222}) for
$W = \inf\limits_{k \ge 1} \frac{\left(1+\varepsilon\right)^k}{k} >
0 $.

\smallskip

\begin{theorem} \label{t3}
Let there exist $ \varepsilon > 0 $ such that
\begin{equation}
\int_0^{\infty}
\left(S\mu(t)-\left(1+\varepsilon\right)\lambda(t)\right)\, dt =
+\infty \label{zkts05}
\end{equation}
\noindent $($large service rates$)$. Then the process $ X (t) $ is
weakly ergodic and has the limiting mean. Moreover, the following
bounds hold:
\begin{equation}
\|{\bf p}(t) - {\bf \pi}(t)\| \le 4\left(1+\varepsilon\right)^k \varepsilon^{-1} e^{-\int\limits_0^t
\frac{\varepsilon}{1+\varepsilon}\left(S\mu(\tau)-(1+\varepsilon)\lambda(\tau)\right)\, d\tau} \label{zkts05a}
\end{equation}
\noindent and
\begin{equation}
\left|E_k(t) - E_0(t)\right| \le
\frac{4\left(1+\varepsilon\right)^k}{\varepsilon W}
e^{-\int\limits_0^t
\frac{\varepsilon}{1+\varepsilon}\left(S\mu(\tau)-(1+\varepsilon)\lambda(\tau)\right)\,
d\tau}, \label{zkts05b}
\end{equation}
\noindent for any initial number of customers $ X (0) = k $, where $
{\bf \pi} (t) $ and $ \phi (t) = E_0(t) $ are the limiting regime
and the limiting mean corresponding to the empty initial queue.
\end{theorem}

{\bf Proof.} Put $d_0=1$, $d_{k+1}=(1+\varepsilon)d_k,\ k \ge 0$,
then instead of (\ref{219}) and (\ref{zkts040}) we have the
following bound of the logarithmic norm:
\begin{eqnarray}
\gamma \left(B(t)\right)_{1D} \le  -\left(\zeta\xi(t) + \frac{\varepsilon}{1+\varepsilon}\left(S\mu(t)-(1+\varepsilon)\lambda(t)\right)\right) \le \\\nonumber
 -\frac{\varepsilon}{1+\varepsilon}\left(S\mu(t)-(1+\varepsilon)\lambda(t)\right) = -\alpha_*(t). \label{zkts06}
\end{eqnarray}
This estimate implies our claim.

\begin{remark} Perturbation bounds for general inhomogeneous $M_t|M_t|S$ queue with
catastrophes can be formulated using the general approach of
\cite{z14} and previous bounds on the rate of convergence.
\end{remark}

\section{Truncations}

Unfortunately, the structure of the infinitesimal matrix of the
process does not provide uniform truncation bounds, as in
\cite{zsks14}. Instead, we can apply another approach to finding
simple and sufficiently sharp truncation bounds, the first such
example was considered in \cite{z11}.

Consider the family of ``truncated'' processes $X_{n}(t)$ on the
state space $E_{n} = \{0,1,\dots,n\}$ with the corresponding reduced
intensity matrix $A_n(t)$. Below we will identify the finite vector
with entries, say,  $ (a_1, a_2, \dots, a_n)^{T}$ and the infinite
vector with the same first $n$ coordinates and the others equal to
zero. In addition, we suppose that
\begin{equation}
e^{-\int_s^t\alpha(u)\, du} \le M e^{-a(t-s)}, \label{4001}
\end{equation}
\noindent for some positive $M,a$ and any $s,t,\ 0 \le s \le t$. Put
$W_n=\inf_{k \ge n} \frac{\sum_{i=n}^{k} d_i}{k}$.

\begin{theorem} \label{t4}
Let the assumptions of Theorem \ref{t1} be fulfilled, and, in
addition,  let (\ref{4001}) hold. Then
 \begin{eqnarray}
\left\|{\bf p} (t) - {\bf p}_n (t)\right\| \le  \frac{8L t}{n W_n }
\left( M jd_{j+1} + \frac{LMd_1}{a}\right) ,
\label{4002}\end{eqnarray}
\begin{eqnarray}
\left|E_{\bf p}(t)- {E}_{\bf p_n}(t)\right| \le  \frac{3L(n+1)t}{n
W_n}\left( M jd_{j+1} + \frac{LMd_1}{a}\right). \label{4003}
\end{eqnarray}\noindent for any $t \ge 0$, and initial condition ${\bf p}(0)={\bf p}_n(0)={\bf e}_j$.
\end{theorem}

{\bf Proof.}  Consider the forward Kolmogorov equation for $X(t)$
and $X_n(t)$ respectively  in the following form:
\begin{equation}
\frac{d\mathbf{p}}{dt}=A_n(t) \mathbf{p}  +\left(A(t) -
A_n(t) \right) \mathbf{p}, \label{cat06}
\end{equation}
\noindent and
\begin{equation}
\frac{d\mathbf{p_n}}{dt}=A_n(t) \mathbf{p_n}.
\label{cat061}
\end{equation}
We have
\begin{equation}
{\bf p}_n (t)= U_n(t,0){\bf p} (0)
\end{equation}
\noindent if ${\bf p} (0) = {\bf p}_n (0)$ and
\begin{equation}
{\bf p} (t)= U_n \left(t,0\right) {\bf p} (0) + \int\limits_0^t U_n \left(t,
\tau\right) \left(A(\tau) - A_n(\tau) \right) {\bf p} (\tau)\,
d\tau. \label{eq3316}
\end{equation}
Then in {\it any} norm  we have
\begin{equation}
\left\|{\bf p} (t) - {\bf p}_n (t)\right\| = \left\|\int\limits_0^t
U_n \left(t,
 \tau\right) \left(A(\tau) - A_n(\tau) \right) {\bf p} (\tau)\, d\tau \right\|.
\label{cat063}
\end{equation}
Consider the Cauchy matrix
\begin{equation}
U_n =\left(
  \begin{array}{ccccccc}
  u_{00}^n & . & . & u_{0n}^n  & 0 & 0 & \cdots \\
u_{10}^n & . & . & u_{1n}^n  & 0 & 0 & \cdots \\
\cdots \\
u_{n0}^n & . & . & u_{nn}^n  & 0 & 0 & \cdots \\
0 & . & . & 0 & 1 & 0 & \cdots \\
0 & . & . & 0 & 0 & 1 & \cdots \\
\cdots
  \end{array}
  \right).
\label{18}
\end{equation}
Then
\begin{eqnarray*}
\left(A -A_n\right) {\bf p} =
\end{eqnarray*}
\begin{scriptsize}
\begin{eqnarray*}
= \left(\sum_{i>n}\xi_i(t) p_i,0,\dots, -\lambda_n(t)p_{n}+ \mu_{n+1}(t)p_{n+1},
\lambda_n(t)p_{n}-(\lambda_{n+1}(t)+\mu_{n+1}(t)+\xi_{n+1}(t))p_{n+1}+\mu_{n+2}(t)p_{n+2}, \dots \right)^T
\end{eqnarray*}
\end{scriptsize}
\noindent and hence
\begin{equation}
U_n\left(A -A_n\right) {\bf p} =\left(\begin{array}{c}
u_{00}^n\sum_{i>n}\xi_i(t) p_i+ u_{0n}^n(-\lambda_n(t)p_{n}+ \mu_{n+1}(t)p_{n+1})\\
u_{10}^n\sum_{i>n}\xi_i(t) p_i+ u_{1n}^n(-\lambda_n(t)p_{n}+ \mu_{n+1}(t)p_{n+1})\\
\vdots \\
u_{n0}^n\sum_{i>n}\xi_i(t) p_i + u_{nn}^n(-\lambda_n(t)p_{n}+ \mu_{n+1}(t)p_{n+1})\\
\lambda_n(t)p_{n}-(\lambda_{n+1}(t)+\mu_{n+1}(t)+\xi_{n+1}(t))p_{n+1}+\mu_{n+2}(t)p_{n+2} \\
\lambda_{n+1}(t)p_{n+1}-(\lambda_{n+2}(t)+\mu_{n+2}(t)+\xi_{n+2}(t))p_{n+2}+\mu_{n+3}(t)p_{n+3}\\ \vdots
\end{array}\right). \label{21}
\end{equation}
\begin{eqnarray}
\|U_n \left(A -A_n \right) {\bf p}\| \le \nonumber \\
 \sum_{k\ge 0}^n\left| u_{k0}^n\sum_{i>n}\xi_i(t) p_i\right|+ \sum_{k\ge 0}^n \left|u_{kn}^n(-\lambda_n(t)p_{n}+ \mu_{n+1}(t)p_{n+1})\right|+ \nonumber \\
 \sum_{k\ge n}\left|\lambda_{k}(t)p_{k}-(\lambda_{k+1}(t)+\mu_{k+1}(t)+\xi_{k+1}(t))p_{k+1}+\mu_{k+2}(t)p_{k+2} \right| \le \nonumber \\
L\sum_{i>n} p_i+\left|\lambda_n(t)p_{n}\right|+ \left|\mu_{n+1}(t)p_{n+1}\right|+ \nonumber \\
2\sum_{k\ge n}\left|\lambda_{k}(t)p_{k}\right|+2 \sum_{k\ge n+1}\left| \mu_{k}(t)p_{k}\right|+ \sum_{k\ge n+1} \left| \xi_{k}(t)p_{k}\right| \le \frac{8L}{n}\sum_{k\ge n}kp_{k},
\end{eqnarray}
\begin{eqnarray}
\|U_n \left(A -A_n \right) {\bf p}\|_{1E} = \nonumber \\
\sum_{k\ge 1}^n k\left|  u_{k0}^n\sum_{i>n}\xi_i(t) p_i\right|+ \sum_{k\ge 1}^n k\left| u_{kn}^n(-\lambda_n(t)p_{n}+ \mu_{n+1}(t)p_{n+1})\right|+ \nonumber \\
 \sum_{k\ge n} (k+1) \left|\lambda_{k}(t)p_{k}-(\lambda_{k+1}(t)+\mu_{k+1}(t)+\xi_{k+1}(t))p_{k+1}+\mu_{k+2}(t)p_{k+2}\right| \le \nonumber \\
  L\sum_{k\ge n}{\left(2k+1\right)p_k}\le \frac{3L(n+1)}{n}  \sum_{k \ge n}kp_{k}.
\end{eqnarray}
On the other hand,
\begin{eqnarray}
\| {\bf p} \|_{1D} \ge d_{n} (p_n + p_{n+1} + \cdots)+ d_{n+1} (p_{n+1} + p_{n+2} + \cdots)+ \cdots =  \nonumber \\
 p_n d_{n}+ p_{n+1}(d_{n}+d_{n+1})+\cdots = \nonumber \\
\frac{d_{n}}{n}np_n+\frac{d_{n}+d_{n+1}}{n+1}p_{n+1}+ \cdots \ge W_n \sum_{k \ge n} {k p_k}.
\end{eqnarray}
Therefore we have
\begin{eqnarray}
\| {\bf p} \|_{1D} \le
\|V(t,0){\bf p}(0) \|_{1D} +
 \int\limits_0^t \| V(t,\tau){\bf f}(\tau)\, d\tau \|_{1D} \le \nonumber \\
e^{-\int\limits_{0}^t \alpha(u)\, du}\|{\bf p}(0) \|_{1D} +  \int\limits_0^t \lambda_0(\tau)
e^{-\int\limits_{\tau}^t \alpha(u) \, du} \, d\tau \le \\
  Me^{-at}\|{\bf p}(0) \|_{1D} + \frac{LMd_1}{a} \le  M jd_{j+1} + \frac{LMd_1}{a} \nonumber
\label{3051}
\end{eqnarray}
\noindent for any ${\bf p}(0)={\bf e}_j$, since $\|{\bf
f}(\tau)\|_{1D} \le d_1 L$, and bounds (\ref{4002}), (\ref{4003})
hold.

\smallskip

\begin{corollary}   \label{C2}
Under the assumptions of Theorem \ref{t2}, let there exist positive
$M,a$ such that $ e^{-\int_s^t \left(\zeta
\xi(u)-\varepsilon\lambda(u)\right)\, du} \le M e^{-a(t-s)},$ for
any $s,t,\ 0 \le s \le t$ $($large catastrophe rates$)$. Then the
following bounds hold:
\begin{equation}
\left\|{\bf p} (t) - {\bf p}_n (t)\right\| \le  \frac{8L t}{n W_n }
\left( M j(1+\varepsilon)^{j+1} + \frac{LM(1+\varepsilon)}{a}\right)
,
 \label{3304}
\end{equation}
\noindent
\begin{equation}
\left|E_{\bf p}(t)- {E}_{\bf p_n}(t)\right| \le  \frac{3L(n+1)t}{n
W_n}\left(Mj(1+\varepsilon)^{j+1} +
\frac{LM(1+\varepsilon)}{a}\right),
 \label{3305}
\end{equation}
\noindent for any $t \ge 0$, and any initial condition ${\bf
p}(0)={\bf p}_n(0)={\bf e}_j$.
\end{corollary}

\smallskip

\begin{corollary}   \label{C3}
Under the assumptions  of Theorem \ref{t3}, let there exist $M,a$
such that
\begin{equation}
e^{-\int_s^t\left(S\mu(u)-\left(1+\varepsilon\right)\lambda(u)\right)\,
du} \le M e^{-a(t-s)}, \label{41}
\end{equation}
\noindent for  any $s,t,\ 0 \le s \le t$ $($large service rates$)$.
Then the following bounds hold:
\begin{equation}
\left\|{\bf p} (t) - {\bf p}_n (t)\right\| \le  \frac{8L t}{n W_n }
\left( M j(1+\varepsilon)^{j+1} + \frac{LM(1+\varepsilon)}{a}\right)
,
 \label{3307}
\end{equation}
\noindent and
\begin{equation}
\left|E_{\bf p}(t)- {E}_{\bf p_n}(t)\right| \le  \frac{3L(n+1)t}{n
W_n}\left( M j(1+\varepsilon)^{j+1} +
\frac{LM(1+\varepsilon)}{a}\right),
 \label{3308}
\end{equation}
\noindent for any $t \ge 0$,  and any initial condition ${\bf
p}(0)={\bf p}_n(0)={\bf e}_j$.
\end{corollary}

\section{Example}

Consider an $M_t|M_t|S$ queue with catastrophes in the case of large
$S$ and periodic intensities. Let   $S=10^{12}$,  $\lambda (t) = 1+
\sin 2\pi t$, $\mu(t) = 3+2\cos 2\pi t$, $\xi_k(t)=\zeta_k
 \xi(t)$, where $\xi (t) = 1- \sin 2\pi t$ and $\zeta_k=1+1/k$.
A similar example without catastrophes  was considered in \cite{z13}
and \cite{zsks14}.

Here we briefly discuss the way for choosing $\{d_i\}$. Firstly, the
monotonicity of this sequence implies the bounds
$$\alpha_{0}\left( t\right) \le \mu(t) +\xi_{1}(t), d_1=1,
$$ hence
$$\alpha_{1}\left( t\right) \le \mu(t) +\xi_{2}(t), d_2=1,
$$
and so on. Therefore, the best possible bound of the ``decay
function'' is
$$\alpha^*\left( t\right) = \mu(t) + \inf_{k\geq 0}
\xi_{k}(t) = \mu(t) +  \xi (t) = 4 +2\cos2\pi t- \sin2\pi t.$$ On
the other hand, such approach yields small values of $W$ and $W_n$.
Therefore, we obtain  bad scores both for the rate of convergence to
the limiting mean and for the error of truncations.

Finally we choose the ``average'' sequence $\{d_i\}$, namely,
putting $d_{k+1}=2^k$ for any $k \ge 0$, we have
$$W=\inf_{i \ge 1} \frac {d_i}{i} =1, \quad g_k=\sum_{n=1}^k d_n \le 2^k,$$
$$
\alpha_{k}\left( t\right) \ge \mu(t) +\xi(t) -
 \lambda(t).$$
Therefore, Theorem \ref{t1} gives us the weak ergodicity of  $X(t)$.
Moreover, if the limiting regime and limiting mean correspond to the
initial condition  $X(0)=0$, then ${\bf p^{**}}(0) = {\bf e}_0$ and
$\phi (0)=0$, and the following bounds hold:
\begin{equation}
\|{\bf p^*}(t) - {\bf p^{**}}(t)\|  \le 2^{k+2}e^{-\int\limits_0^t
\left(3 + 2\cos 2 \pi \tau - 2 \sin 2\pi\tau\right)\, d\tau} \le
2^{k+4}e^{-3t}, \label{601}
\end{equation}
\begin{equation}
\left|\phi(t) - E_k(t)\right| \le 2^{k+2}e^{-\int\limits_0^t \left(3
+ 2\cos 2 \pi \tau - 2 \sin 2\pi \tau \right)\, d\tau} \le
2^{k+4}e^{-3t}, \label{602}
\end{equation}
\noindent for any $ t \ge s \ge 0$ and any initial condition
$X(0)=k$.

Consider the error of truncations. We have $M \le 4$, $a=3$, $L
\approx 5\cdot 10^{12}$, and $W_n=\inf_{k \ge n}
\frac{\sum_{i=n}^{k} d_i}{k} = \frac{2^{n-1}}{n}$. Hence, the
following bounds follow from Theorem \ref{t4}:
 \begin{eqnarray}
\left\|{\bf p} (t) - {\bf p}_n (t)\right\| \le  \frac{t \cdot
10^{13}}{2^{n-3}} \left( k 2^{k+2} + 10^{14}\right),
\label{603}\end{eqnarray}
\begin{eqnarray}
\left|E_{\bf p}(t)- {E}_{\bf p_n}(t)\right| \le  \frac{t(n+1) \cdot
10^{14}}{2^{n-1}} \left(  k 2^{k+2} + 10^{14}\right),
\label{604}
\end{eqnarray}\noindent for any $t \ge 0$, and any initial condition $X(0)=k$.

Therefore, we can choose $n=120, \quad t \in [6,7]$ and find the
limiting characteristics with error less then $10^{-6}$.
\begin{figure}[h!]
\begin{center}
\includegraphics[width=8cm]{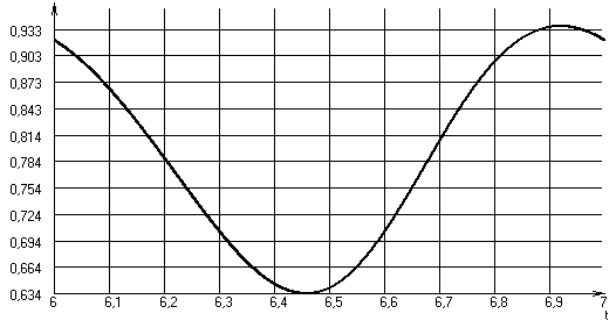}
\end{center}
\caption{Approximation of the limiting probability of empty queue
$\Pr\{X(t) = 0|X(0)=0\}$ on $[6,7]$.}
\end{figure}
\begin{figure}[h!]
\begin{center}
\includegraphics[width=8cm]{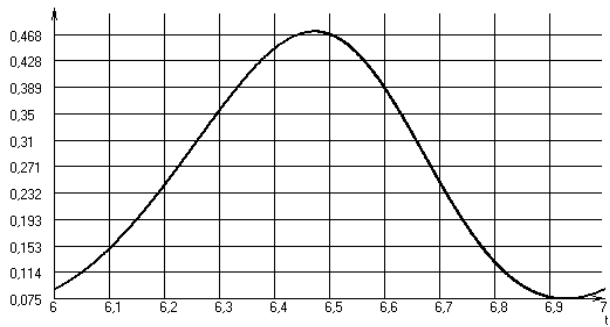}
\end{center}
\caption{Approximation of the limiting mean $E(t,0)$ on $[6,7]$.}
\end{figure}

\smallskip

{\bf Acknowledgement.} This work was supported by Russian Scientific
Foundation  (Grant No. 14-11-00397).

\end{document}